\documentclass{article} 
\usepackage{amsmath,amsfonts,bm}
\usepackage{hyperref}
\usepackage{url}
\usepackage{float}

\usepackage{algorithm}
\usepackage{algpseudocode}
\usepackage{cleveref}
\usepackage{booktabs}
\usepackage{graphicx}
\usepackage{multirow}

\title{A Short Review of Estimators for the GLM predictive of Laplace Bayesian Neural Networks}
\author{Romie Banerjee}
\date{}

\usepackage{fancyhdr}  
\usepackage{lipsum}    

\usepackage{fancyhdr}
\begin{document}
\maketitle
\begin{abstract}
This short review examines the primary approaches for estimating the predictive distribution of Laplace-approximated Bayesian neural networks, with particular focus on the Generalized Linear Model (GLM) formulation. We survey the landscape of estimation strategies, from exact GLM computations requiring full Jacobian evaluations to Monte Carlo approximations that trade computational cost for statistical efficiency. The review covers the theoretical foundations of the Laplace approximation, the Kronecker-factored approximate curvature (KFAC) method for scalable posterior inference, and the various predictive estimation techniques developed in the literature. We provide a unified presentation that clarifies the relationships between methods and highlights their respective computational and statistical trade-offs.
\end{abstract}

\section{Introduction}

\subsection{Motivation}

Deep neural networks have achieved remarkable success across numerous domains, yet their deployment in safety-critical applications demands reliable uncertainty quantification. Bayesian neural networks (BNNs) offer a principled framework for obtaining predictive distributions, but exact Bayesian inference remains computationally intractable for modern architectures. This has motivated a rich literature on approximation strategies, among which the Laplace approximation stands out for its simplicity and post-hoc applicability to pre-trained models.

Given a pre-trained model, one can \textit{bayesianify} it via the Laplace approximation. This additional structure transforms the model from a point-estimate predictor into one that outputs full predictive distributions. The process involves two key components:

\begin{enumerate}
    \item \textbf{Posterior Approximation:} The posterior over network weights is approximated by a Gaussian distribution using the Laplace approximation. This reduces to calculating the Hessian of the loss function at the maximum a posteriori (MAP) estimate.
    \item \textbf{Predictive Estimation:} The predictive distribution, which requires integrating over the posterior, admits no simple closed-form solution due to the nonlinearity of the feed-forward function. Various estimation strategies have been developed to address this challenge.
\end{enumerate}

\subsection{The Laplace Approximation and Its Variants}

The Laplace approximation for neural networks was pioneered by MacKay \cite{mackay}, who recognized that the local curvature of the loss landscape could provide a Gaussian approximation to the posterior distribution. The Hessian computation, however, presents significant computational challenges for modern networks with millions of parameters.

Several approximations have been proposed to make the approach scalable. The Generalized-Gauss-Newton (GGN) approximation replaces the Hessian with a product of first-order derivatives, corresponding to the empirical Fisher information matrix. The Kronecker-factored approximate curvature (KFAC) method \cite{kfac, kfac-cnn} further exploits the structure of neural network gradients, representing layer-wise Fisher blocks as Kronecker products of smaller matrices. This two-level approximation---GGN followed by KFAC---has become the standard for practical Laplace BNNs.

\subsection{Predictive Estimation Strategies}

The predictive distribution of a Laplace BNN is defined as the push-forward of the Gaussian posterior through the network function. Since the network is nonlinear, this distribution is generally non-Gaussian and intractable. The literature has developed two primary estimation strategies:

\begin{enumerate}
    \item \textbf{Monte Carlo Integration:} Samples are drawn from the posterior distribution and propagated through the network. The empirical distribution of the resulting outputs provides an estimator for the true Laplace-BNN predictive. This approach, while straightforward, can be computationally expensive and may not preserve the MAP prediction as the mean of the predictive distribution.
    
    \item \textbf{Generalized Linear Model (GLM):} The network function is linearized around the MAP estimate, yielding a tractable Gaussian predictive distribution. This approach, introduced in the context of Laplace BNNs by Immer et al. \cite{glm}, provides a closed-form solution and preserves the MAP prediction as the predictive mean.
\end{enumerate}

\subsection{Estimating the GLM Predictive}

Given an input $x$, a pre-trained network with parameters $\theta_*$, and the network function $f_{\theta_*}(-)$, the GLM predictive is given by

\begin{equation}
P(y|x) = \mathcal{N}\left(f_{\theta_*}(x),\, J(x)\,\Sigma\,J(x)^T\right),
\end{equation}

where $\Sigma$ is the covariance of the Laplace-approximated posterior and $J(x)$ is the Jacobian of $f_\theta(x)$ with respect to $\theta$ evaluated at $\theta_*$. 

The covariance $\Sigma$ is computed offline and remains fixed for all inputs, whereas the Jacobian $J(x)$ must be computed online for each new input. This presents a significant computational bottleneck, particularly when $f_\theta$ is high-dimensional or when the output space is large.

\subsection{Review Scope and Contributions}

This short review synthesizes the literature on estimators for the GLM predictive of Laplace BNNs. We examine:

\begin{enumerate}
    \item The exact GLM formulation and its computational requirements
    \item The KFAC approximation and its impact on predictive estimation
    \item Monte Carlo sampling strategies for approximating the GLM predictive
    \item Computational trade-offs between exact and approximate methods
    \item Practical implementation considerations and convergence properties
\end{enumerate}

\section{Bayesian Neural Networks}

A Bayesian neural network (BNN) consists of a (differentiable) map $f:X \times \Theta \to Y$, where the spaces $X$, $\Theta$ and $Y$ and $Y$ are respectively the input, model weights and output spaces. The spaces are isomorphic as affine spaces to $X \simeq \mathbb{R}^n$, $\Theta \simeq \mathbb{R}^w$ and $Y \simeq \mathbb{R}^d$. The mapping $(x,\theta) \mapsto f(x,\theta)$ expresses the evaluation of a model $\theta$ on the input $x$. (Note: In this formulation, the aleatoric predictive is suppressed by considering on point-valued functions $f$.). In this setup, a regular feed-forward network is a BNN with a fixed choice of weights, for example $\theta = \theta_*$. We can denote the resulting network as a function $f_{\theta_*}:X \to Y$, $f_{\theta_*}(x) = f(x, \theta_*)$. 
The weights space $\Theta$ is equipped with a probability measure, the posterior distribution (usually intractable) $P(\theta | D)$. The mode of the posterior, 

\begin{equation}
\theta_{\text{MAP}} = \text{argmin}_{\theta \in W} \log P(D|\theta)  + \log P(\theta)
\end{equation}
is usually estimated using empirical risk minimization. Given.e..g., an i.i.d. classification data-set $D:= {(x_n, y_n) \in X \times Y}_{n=1}^N$, the weights $\theta\in \Theta$ are trained to minimize the regularized empirical risk $\mathcal{L}(D;\theta) = \log P(D|\theta)  + \log P(\theta)$.

The predictive distribution $P(y|x)$ on the output space $Y$ is the push-forward probability measure of $P(\theta|D)$ under the map $f(x, -):\Theta \to Y$.  
\begin{equation}
P(y|x) = f(x,-)_* P(\theta|D) 
\end{equation}

The (bayesian) prediction $\hat{f}(x)$ is the mean of $P(y|x)$, i.e. 
\begin{equation}
\hat{f}(x) =\mathbb{E}_{y \sim P(y|x)}[y] =  \mathbb{E}_{\theta \sim P(\theta|D)}\left[ f(x, \theta) \right]
\end{equation}

The predictive uncertainty is the (co)variance of the predictive distribution $P(y|x)$, i.e., 

\begin{equation}
\text{unc}(\hat{f}(x)) = \text{cov}_{y \sim P(y|x)}[y] =  \text{cov}_{\theta \sim P(\theta|D)} \left[ f(x, \theta) \right]
\end{equation}

\section{Monte Carlo Integration for BNNs}
The Monte-Carlo estimate of predictive distribution is obtained by first sampling weights $\theta_1, \cdots, \theta_k \sim P(\theta|D)$. This produces samples from the predictive distribution $\{ f(x, \theta_i) \}_{i=1}^k \sim P(y|x)$. Compute the sample mean to estimate the Bayesian prediction 
\begin{equation}
\overline{\hat{f}(x)} = \frac{1}{k} \sum_{i=1}^{k} f(x, \theta_i) \in \mathbb{R}^d
\end{equation}
Compute the sample covariance to estimate the predictive predictive 

\begin{equation}
\overline{\text{unc}(\hat{f}(x))} = \frac{1}{k-1} \sum_{i=1}^k f(x, \theta_i) f(x, \theta_i)^T \in \mathbb{R}^{d\times d}
\end{equation}

One successful method following this principle is MCDropout \cite{gal2015bayesian, gal2016dropout}, where the posterior is estimated by variational inference on Bernoulli distributions. 
One disadvantage of the MC sampling method is that it is not post-hoc as the estimated $\overline{\text{unc}(\hat{f}(x))}$ is an predictive for the Bayesian prediction $\overline{\hat{f}(x)}$ which in general does not equal to the original MAP prediction $f(x, \theta_{\text{MAP}})$.

\section{Laplace Approximation}
(\cite{mackay}) Laplace approximation replaces the bayesian posterior with a normal distribution with $\theta_{\text{MAP}}$ as mean, $P(\theta | D) \sim_{LA} \mathcal{N}(\theta_{\text{MAP}}, \Sigma)$, where the covariance  $\Sigma\in  \mathbb{R}^{w\times w}$ is the inverse of the Hessian 
\begin{equation}
\Sigma := \left(\left.\frac{\partial^2}{\partial\theta^2} \mathcal{L}(D;\theta)\right\vert_{\theta_{\text{MAP}}}\right)^{-1}
\end{equation}

This is second-order derivative calculation is often replaced in practice by the products of first-order derivatives , i.e. Generalized-Gauss-Newton approximation or equivalently the Fisher information matrix; the covariance of the gradients of the loss function at the empirical risk minimizer $\theta_{\text{MAP}}$.

\begin{align*}
I_{\theta}(\theta_{\text{MAP}}) = \text{Cov}_{(x,y) \sim D} \left(\left.\frac{\partial}{\partial\theta} \mathcal{L}((x,y); \theta)\right\vert_{\theta_{\text{MAP}}}\right) 
\end{align*}

\begin{equation}
\boxed{\Sigma \sim I_{\theta}(\theta_{\text{MAP}})^{-1}}
\end{equation}

\section{Generalized Linear Model for Laplace-BNNs}

(\cite{glm}) Define the generalized linear model by linearizing the BNN $f:X \times \Theta \to Y$ at $\theta = \theta_{\text{MAP}}$ (via. Taylor series)

\begin{equation}
f_{\text{lin}}(x, \theta) = f(x, \theta_{\text{MAP}}) + \left[\left.\frac{\partial}{\partial\theta} f(x, \theta)\right\vert_{\theta_{\text{MAP}}}\right]  \left(\theta - \theta_{\text{MAP}}\right)
\end{equation}

Denote the jacobian matrix  $\left.\frac{\partial}{\partial\theta} f(x, \theta)\right\vert_{\theta_{\text{MAP}}}$ by $J(x) \in \mathbb{R}^{d \times w}$. Replace the posterior distribution $P(\theta|D)$ by the Laplace approximation $\mathcal{N}(\theta_{\text{MAP}}, \Sigma)$. The GLM predictive distribution $P_{\text{lin}}(y|x)$ on the output space $Y$ is defined to be the push-forward measure of the Laplace posterior $N(\theta_{\text{MAP}}, \Sigma)$ on the weight space $\Theta$ under the map $f_{\text{lin}}(x, -): \Theta \to Y$. Since $f_{\text{lin}}(x, \theta)$ is affine on $\theta$, the predictive distribution $P_{\text{lin}}(y|x)$ is normal.

The mean and covariance of the normal distribution $P_{\text{lin}}(y|x)$ are as follows

\begin{align*}
\mathbb{E}_{y \sim P_{\text{lin}}(y|x)}[y] &= \mathbb{E}_{\theta \sim \mathcal{N}(\theta_{\text{MAP}}, \Sigma)} \left[ f_{\text{lin}}(x, \theta) \right] \\ 
&= f(x, \theta_{\text{MAP}}) + J(x) *\mathbb{E}_{\theta \sim \mathcal{N}(\theta_{\text{MAP}}, \Sigma)}\left[ (\theta - \theta_{\text{MAP}})\right] \\
&=  f(x, \theta_{\text{MAP}})
\end{align*}

\begin{align*}
\text{cov}_{y \sim P_{\text{lin}}(y|x)}[y] &= \text{cov}_{\theta \sim \mathcal{N}(\theta_{\text{MAP}}, \Sigma)} \left[ f_{\text{lin}}(x, \theta) \right]\\
&= \text{cov}_{\theta \sim \mathcal{N}(\theta_{\text{MAP}}, \Sigma)} \left[J(x) * (\theta - \theta_{\text{MAP}})\right]\\
&= \text{cov}_{\eta \sim \mathcal{N}(0, \Sigma)} \left[J(x) * \eta\right]\\
&= J(x) *\Sigma*J(x) ^T\\
\end{align*}

This property makes GLM suitable for post-hoc predictive estimation as the Bayesian prediction $\hat{f}_{\text{lin}}(x)$ which is simply the mean of the predictive distribution, agrees with the original prediction from the model $\theta_{\text{MAP}}$. 

\begin{equation}
\hat{f}_{\text{lin}}(x)= f(x, \theta_{\text{MAP}})
\end{equation}

The predictive at this point further has a simple analytical formula, 
\begin{equation}\label{glm_unc}
\boxed{\text{unc}(\hat{f}_{\text{lin}}(x)) = J(x)*\Sigma*J(x)^T}
\end{equation}

\section{Laplace BNNs with KFAC posterior}

(\cite{kfac, kfac-cnn}) The FI matrix $I = I_{\theta}(\theta_{\text{MAP}})$ has a very large size $w \times w$. For practical applications various approximations of $I$ are available. The KFAC method employs two levels of approximation 

\begin{enumerate}
\item Treating each layer on the neural network separately, ignoring cross-layer terms, expressing $I \in \mathbb{R}^{w \times w}$ in block diagonal form with diagonal blocks $I_{(l)} \in \mathbb{R}^{w_l \times w_l},$ where $w_l$ is the size of the $l$-th layer and, $\sum_{l=1}^L w_l = w$. This reduces the size of $I$ from $w^2$ to $\sum w_l^2$. 
\item For each layer, the gradients $\nabla_l$ in this layer, is expressible as a Kronecker product $\nabla_l = a_{l-1}\otimes g_l$, where $a_{l-1}$ is the incoming activation from layer $l-1$ and $g_l$ is the otgoing gradient of layer $l$. The Fisher information block for layer $l$ is then $I_{(l)} = \mathbb{E}\left[\nabla_l*\nabla_l^T\right] = \mathbb{E} \left[(a_{l-1}\otimes g_l)*(a_{l-1}\otimes g_l)^T\right]$. The KFAC method makes the approximation $I_{(l)} \approx \mathbb{E}\left[a_{l-1}\otimes a_{l-1}^T \right] \otimes \mathbb{E}\left[g_l\otimes g_l^T \right]=: Q_{(l)}\otimes H_{(l)}$. Here $Q_{(l)} \in \mathbb{R}^{l_{\text{in}} \times l_{\text{in}}}$ and $H_{(l)} \in \mathbb{R}^{l_{\text{out}} \times l_{\text{out}}}$ and the $l$-th layer is a map $\mathbb{R}^{l_{\text{in}}} \to \mathbb{R}^{l_{\text{out}}}$. This reduces the size from $\sum w_l^2 = \sum (l_{\text{in}}l_{\text{out}})^2$ to $\sum (l_{\text{in}}^2 + l_{\text{out}}^2)$.
\item The kronecker factored block-diagonal form $I_{\text{KFAC}} = \text{diag} \left [Q_{(1)}\otimes H_{(1)}, \cdots, Q_{(L)} \otimes H_{(L)} \right]$.
\item The KFAC posterior covariance, $\Sigma_{\text{KFAC}}$ has the form
\begin{equation}\label{sigma_kfac}
\begin{pmatrix}
Q_{(1)}^{-1}\otimes H_{(1)}^{-1}&  &\\
& \ddots & \\
&& Q_{(L)}^{-1} \otimes H_{(L)}^{-1}
\end{pmatrix}
\end{equation}
\end{enumerate}

\section{The GLM predictive for Laplace-KFAC posterior}

In order to use this during inference one needs access two matrices,  

\begin{itemize}
\item the covariance matrix $\Sigma$ or the inverse of the Fisher information matrix $I_{\theta}(\theta_{\text{MAP}})$ which  is computed only once offline and used repeatedly during inference (independent of the input $x$)
\item the jacobian $J(x)$, which must be computed for every input $x$. 
\end{itemize}

Using the given co-ordinatization $\theta_1, \ldots, \theta_w$ of $\Theta$, the entries of the jacobian matrix $J(x)$ are (\ref{matrix}), where $\frac{\partial f_i}{\partial \theta_j}$ is short for $\left.\frac{\partial}{\partial\theta_j}f_i(x, \theta) \right\vert_{\theta = \theta_{\text{MAP}}}$, and $(f_1, \ldots, f_d)$ are the scalar components of $f$.  The rows are the gradients of the scalar valued functions $f_i, 1\leq i \leq d$.

\begin{equation}\label{matrix}
J(x) = 
\begin{pmatrix}
- & \left[\nabla f_1\right ]^T & - \\
 & \vdots &  \\
- & \left[ \nabla f_d \right]^T &  -
\end{pmatrix}
\end{equation}

During inference time these rows can be computed using back-propagation. This would require $d$ back-propagation steps, one for each scalar value of $f$. 

When the covariance is expressed in Kronecker factored form the GLM uncertianty expression must be 

\begin{equation}\label{glm_unc_kfac}
\text{unc}(\hat{f}_{\text{lin}}(x)) = J(x)_{\text{KF}}*\Sigma_{\text{KFAC}}*J(x)_{\text{KF}}^T
\end{equation}

where $J(x)_{\text{KF}}$ is the jacobian in a kronecker-factored form.

When $d=1$, i.e. $f(x, \theta) \to \mathbb{R}$ is scalar valued, consider the jacobian $\left.\frac{\partial}{\partial\theta} f(x,\theta)\right\vert_{\theta_{\text{MAP}}} \in \mathbb{R}^{1\times w}$. For every layer $l \in [1,\cdots, L]$, the layerwise derivative has a kronecker decomposition $\left.\frac{\partial}{\partial\theta_l} f(x,\theta)\right\vert_{\theta_{\text{MAP}}} = a_{l-1}\otimes g_l$, where $a_{l-1} \in \mathbb{R}^{1\times l_{\text{in}}}$ is the incoming activation and $g_l \in \mathbb{R}^{1 \times l_{\text{out}}}$ the outgoing gradient for layer $l$. This gives a decomposition of the jacobian into blocks.

\begin{align*}
J(x) &= \left( \left.\frac{\partial}{\partial\theta_1} f(x,\theta)\right\vert_{\theta_{\text{MAP}}}, \cdots, \left.\frac{\partial}{\partial\theta_L} f(x,\theta)\right\vert_{\theta_{\text{MAP}}} \right)   \\
&\in \mathbb{R}^{d\times (w_1 + \cdots + w_L)}\\
&= 
\begin{pmatrix}
\vert && \vert \\
a_0^j\otimes g_1^j & \cdots &a_{L-1}^j\otimes g_L^j  \\
\vert && \vert 
\end{pmatrix}
\end{align*}

where the superscript $j \in (1, \cdots, d)$. Finally the GLM-predictive calculation,

\begin{equation*}
J(x)
* 
\begin{pmatrix}
Q_{(1)}^{-1}\otimes H_{(1)}^{-1}&& \\
& \ddots & \\
&& Q_{(L)}^{-1} \otimes H_{(L)}^{-1}
\end{pmatrix}
*
J(x)^T
\end{equation*}

\begin{equation}
=\text{Diag}_{j=1}^d\left[\sum_{l=1}^L \left( a_{l-1}^j*Q_{(l)}^{-1}*{a_{l-1}^j}^T \right) \otimes \left( g_l^j * H_{(l)}^{-1} * {g_{l-1}^j}^T\right)\right]
\end{equation}

\section{A monte-carlo estimation of the GLM predictive}

A Monte-carlo estimator for the GLM predictive is obtained through sampling from the Gaussian posterior.
The jacobian $J(x)$ gets replaced by a low-rank and base-changed version $A(x)$ and the predictive covariance is estimated as $$A(x)*A(x)^T \approx J(x)*\Sigma*J(x)^T$$

The columns of $A(x)$ are the \textit{directional derivatives} of $f(x,\theta)$ along the monte-carlo sample directions, and are approximated numerically. 
Unlike the rows of $J(x)$ which are gradients of the same functions (along standard basis dimensions of the model weight space) and are calculated with \textit{symbolic differentiation}.
This makes computation less reliable on repeated \textit{backprop} calls (for a single $x$) and can be possibly parallelized (due to the fixed sample size). 

The GLM uncterianty term \ref{glm_unc} can be estimated by computing sample covariance,
\begin{align*}
J(x)*\Sigma * J(x)^T &= \text{cov}_{\eta \in \mathcal{N}(0, \Sigma)}\left[ J(x) * \eta \right ] \\
&\approx \frac{1}{N-1} \sum_{i=1}^N \left[J(x)*\eta_i \right] \left[J(x)*\eta_i \right]^T
\end{align*}
where $\left( \eta_1, \cdots, \eta_n \right) \sim \mathcal{N}(0,\Sigma)$ are i.i.d samples. The noise vectors $\eta$ can be interpreted as members of the tangent space $T_{\theta_{\text{MAP}}}(\Theta)$, and the terms in the summand, $J(x)*\eta$ are simply the directional derivatives $D_{\vec{\eta}} f = \lim_{h \to 0} \left(f(\theta_{\text{MAP}} + h\vec{\eta}) - f(\theta_{\text{MAP}})\right)/h$.

The MC-GLM method to estimate predictive is the following:

\begin{enumerate}
\item Sample $(\vec{\eta}_1, \ldots, \vec{\eta}_k) \sim \mathcal{N}(0, \Sigma)^{\times k}$, $k \ll w$. Define,
\item
\begin{equation*}
\overline{A}(x) := 
\begin{pmatrix}
| & & | \\
\bar{D}_{\vec{\eta}_1}f & \cdots & \bar{D}_{\vec{\eta}_k}f \\
| &  & |
\end{pmatrix}
\in \mathbb{R}^{d \times k}
\end{equation*}
\item The approximate directional derivative $\bar{D}$ is calculated using finite differences.
\item $\text{mcglm-unc}(f(x, \theta_{\text{MAP}})):= \overline{A}(x)*\overline{A}(x)^T$
\end{enumerate}

The Cholesky decomposition $\Sigma = B*B^T$ lets us write the predictive $J(x)*\Sigma*J(x)^T$ as $A(x)*A(x)^T$, where $A(x) = J(x)*B$. The matrix $A(x)$ is simply the jacobian $J(x)$ base-changed to the a new basis formed out of the columns of $B$. 
From this point of view the idea of MC-GLM is to base-change $J(x)$ using a different and low-rank basis for $T_{\theta_{\text{MAP}}}(W)$, obtained by sampling. 

\section{Sampling from the KFAC posterior}

In order to generate a sample from a $n$-variable normal $\vec{z} \sim \mathcal{N}(\vec{\mu}, \Sigma)$ one has to write it as $\vec{z} = \vec{\mu} + B\vec{x}$, where $B$ is a $n\times n$ matrix such that $BB^T = \Sigma$ and $\vec{x} \sim \mathcal{N}(\vec{0}, I_n)$. The matrix $B$ can be found using Cholesky decomposition and is guaranteed when $\Sigma$ is a covariance. (Notation: $B = \text{ch}(\Sigma)$.)

When the covariance is in KFAC form as in \ref{sigma_kfac}, we observe that the cholesky decomposition commutes with block-diagonal, inverse and kronecker product. Consequenty the Cholesky is computed as a block diagonal matrix with the blocks $\text{ch}(Q_{(l)})^{-1}\otimes \text{ch}(H_{(l)})^{-1}$ along the diagonal. 
For every layer $l$, let $\mathcal{N}(\vec{\mu}_l, Q_{(l)}^{-1}\otimes H_{(l)}^{-1})$ be the restriction of $\mathcal{N}(\vec{\mu}, \Sigma_{\text{KFAC}})$ restricted to weights of $f_{\theta_l}$.The layer-wise sample obtained as 
\begin{align*}
\vec{z}_l &=  \left( \text{ch}(Q_{(l)})^{-1}\otimes \text{ch}(H_{(l)})^{-1}\right )* \vec{x}_l \\
&= \text{ch}(Q_{(l)})^{-1}*\vec{x}_l\text{.reshape}(l_{\text{in}}, l_{\text{out}})*\text{ch}(H_{(l)})^{-1}
\end{align*}
where  $\vec{x}_l \sim \mathcal{N}(\vec{\mu_l}, I_{w_l})$.

\section{Conclusion}

This review has surveyed the landscape of estimators for the GLM predictive of Laplace Bayesian neural networks. We have examined the theoretical foundations, computational considerations, and practical trade-offs of the primary approaches. The choice of estimator ultimately depends on the specific application requirements, balancing computational budget, desired accuracy, and the need for closed-form solutions versus empirical flexibility.

Future directions include the development of adaptive sampling strategies that allocate computational resources based on input uncertainty, the integration of these methods with modern hardware accelerators, and the exploration of hybrid approaches that combine the strengths of multiple estimation strategies.

\end{document}